\documentstyle[fleqn,12pt]{article}
\begin{document}
\pagenumbering{arabic}
\setcounter{page}{1}
\pagestyle{plain}
\baselineskip=16pt

\thispagestyle{empty}
\rightline{MSUMB 95-04, September 1995}

\begin{center}
{\bf Two parameter deformation of 
     grassmann matrix group \\and supergroup }
\end{center}

\vspace{1cm}
\noindent
Salih Celik\footnote{New E-mail address: sacelik@yildiz.edu.tr}

\noindent
Mimar Sinan University, Department of Mathematics, \\
80690 Besiktas, Istanbul.

\vspace{1cm}
\noindent
{\bf Abstract}

\noindent
The two parameter quantum deformation of 2x2 Grassmann matrices, Gr(2), 
and supermatrices, Gr$(1\vert 1)$, are presented. Gr(2) whose matrix elements 
are all Grassmannian variables is called  the superdual of the genel linear 
group GL(2), and Gr$(1\vert 1)$ whose diagonal matrix elements are 
Grassmannian variables is called the superdual of the supergroup 
GL$(1\vert 1)$ whose nondiagonal elements are Grassmannian. Noncentral dual 
superdeterminant for Grassmann  supermatrices belonging to 
Gr$_{p,q}(1\vert 1)$ is constructed. As with the 2x2 quantum matrices, 
the relations satisfied by the matrix elements of the Grassmann matrices are 
expressed in terms of an $\widehat{R}$-matrix. The properties of $n$th power 
of a Grassmann supermatrix are given as an Appendix. 

\vfill\eject
\noindent
{\bf I. INTRODUCTION }

Quantum groups are a generalization of the concept of groups. 
More precisely, a quantum group is a deformation of a group 
that, for particular values of the deformation parameter, 
coincides with the group. The theory and applications of 
quantum groups have attracted a lot of attention among 
mathematicians and physicists. The main physical motivation for quantum groups 
is that when nonlinear physical systems which are classically completely 
integrable are quantized, th classical symmetry group should be replaced 
by the corresponding quantum group.$^1$ On the other hand, most of the 
difficulties involving the divergences of quantum field theories which lie at 
the heart of all interactions require supersymmetry and thus the introduction 
of supergroups.$^2$ The algebraic structure underlying quantum groups  
extends the theory of the supergroups.$^{3-8}$ In Ref. 9, the $q$-analog, 
Gr$_q(1\vert 1)$, of the dual supermatrices Gr$(1\vert 1)$ is presented. 
Gr$_q(1|1)$ is the superdual of the quantum group GL$_q(1|1)$ and the 
properties of the quantum dual supermatrices are discussed. 
In this paper we present a two parameter deformation, Gr$_{p,q}(2)$ and 
Gr$_{p,q}(1\vert 1)$, of the Grassmann matrices and supermatrices, 
respectively and give an $\widehat{R}$-matrix for this deformation. 

We will say that Grassmann matrices are the dual matrices in 
GL(2) and Grassmann supermatrices are the dual supermatrices 
in GL$(1|1)$. To study the two parameter extension of the 
Grassmann matrices and supermatrices, we follow the approach 
of Manin$^3$ in Sec. II. In the following section we get an 
$\widehat{R}$-matrix which gives the relations between the 
matrix elements of a dual matrix in GL$(1|1)$. The properties 
of the $n$-th power of a dual supermatrix which are more compact 
than the single deformation parameter case are given in Appendix. 

\noindent
{\bf II. THE QUANTUM GRASSMANN MATRIX GROUP Gr$_{p,q}(2)$ }

Before discussing the two parameter deformation of the 
dual matrices in the general linear group GL(2), we 
give some notations and useful formulas. 

\noindent
{\bf A. Notations}

Consider 2x2 matrices with Grassmannian entries. We will say 
that such matrices form the Grassmann matrix group and denote 
it by Gr(2). Explicitly, a Grassmannian 2x2 matrix 
$\widehat{A}$ is of the form 
$$\widehat{A} = \left ( \matrix{ \alpha & \beta \cr
                                 \gamma & \delta \cr}
\right) \eqno(2.1) $$
where all entries are Grassmannian. 

Since the matrix elements of $\widehat{A}$ are all Grassmannian, 
for the conventional tensor products 
$$ \widehat{A}_1 = \widehat{A} \otimes I ~~\mbox{and}~~ 
\widehat{A}_2 = I \otimes \widehat{A} \eqno(2.2) $$
one can write (no-grading)
$$(\widehat{A}_1)^{ij} \hspace*{0.0001cm}_{{kl}} = 
 \widehat{A}^i \hspace*{0.0001cm}_{k} \delta^j\hspace*{0.0001cm}_l, 
    \eqno(2.3a)$$
$$(\widehat{A}_2)^{ij} \hspace*{0.0001cm}_{{kl}} = 
 \delta^i\hspace*{0.0001cm}_k \widehat{A}^j \hspace*{0.0001cm}_{l} 
 \eqno(2.3b)$$
where $\delta$ denotes the Kronecker delta. 

\noindent
{\bf B. Two parameter deformation of Gr(2)}
 
The one parameter deformation of Grassmann matrices was given 
by Corrigan et al.$^4$ In this section, we will give a two 
parameter deformation Grassmann matrices, i.e., of $Gr(2)$. 
Let $R_p[2|0]$ be a quantum vector space which is two dimensional. The 
coordinates of a vector $V = (x, y)^T \in R_p[2|0]$ satisfy the bilinear 
product relation 
$$xy - p yx = 0. \eqno(2.4)$$
We consider a dual quantum vector space $R_q[0|2]$, the generators 
of which are Grassmannian. The coordinates of a (dual) vector 
$\widehat{V} = (\xi, \eta)^T \in R_q[0|2]$ satisfy the relations 
$$ \xi^2 = 0 = \eta^2, ~\eta \xi + q \xi \eta = 0 \eqno(2.5) $$
as introduced in Ref. 3. 

Now we want to define a two parameter deformation of the 
algebra of functions on the Grassmann matrix group Gr(2) as 
an associative algebra with unit, generated by the generators 
$\alpha$, $\beta$, $\gamma$, and $\delta$. For this, we 
considering linear transformations $\widehat{A}$ with the 
following properties: 
$$ \widehat{A} : R_p[2|0] \longrightarrow R_q[0|2], \eqno(2.6\mbox{a})$$
$$ \widehat{A} : R_q[0|2] \longrightarrow R_p[2|0]. \eqno(2.6\mbox{b}) $$
We assume that the matrix elements of $\widehat{A}$ commute 
with the coordinates of $R_p[2|0]$ and anti-commute with 
the coordinates of $R_q[0|2]$. Then the endomorphisms in (2.6) 
impose the following $(p,q)$-anti-commutation relations among the 
matrix elements of $\widehat{A}$: 
$$ \alpha \beta + p^{-1} \beta \alpha = 0, \qquad 
   \alpha \gamma + q^{-1} \gamma \alpha = 0, $$
$$ \gamma \delta + p^{-1} \delta \gamma = 0, \qquad 
   \beta \delta + q^{-1} \delta \beta = 0, \eqno(2.7)$$
$$ \alpha \delta + \delta \alpha = 0, \qquad 
   \alpha^2 = \beta^2 = \gamma^2 = \delta^2 = 0, $$
$$ \beta \gamma + p q^{-1} \gamma \beta = (p - q^{-1}) \delta \alpha $$
where $p$ and $q$ are non-zero complex numbers and $pq \pm 1 \neq 0$. 

Since the entries of $\widehat{A}$ are all Grassmannian, a proper 
inverse cannot exist. However, the left and right inverses of 
$\widehat{A}$ can be constructed. Let 
$$\Delta_L = \beta \gamma + q^{-1} \delta \alpha, \eqno(2.8a)$$
$$ \Delta_R = \gamma \beta - p^{-1} \alpha \delta. \eqno(2.8b) $$
Then at least the formally, the left and right inverses of $\widehat{A}$ 
become 
$$\widehat{A}_L^{-1} = \left ( \matrix{ q^{-1}\delta & \beta \cr
                                    - pq^{-1} \gamma & -p \alpha \cr}
\right), \eqno(2.9) $$
$$\widehat{A}_R^{-1} = \left ( \matrix{-q \delta & \beta \cr
                                 - qp^{-1} \gamma & p^{-1} \alpha \cr}
\right). \eqno(2.10) $$
Indeed, it is easy to show that
$$ \widehat{A}_L^{-1} \widehat{A} = \Delta_L I,  \eqno(2.11a)$$
$$ \widehat{A} \widehat{A}_R^{-1} = \Delta_R I \eqno(2.11b) $$
where I is the 2x2 unit matrix. In this case, $\Delta_L$ may be 
considered as a left quantum (dual) determinant and $\Delta_R$ as 
a right quantum (dual) determinant. Note that, one can write 
$$ \Delta_L \widehat{A}_R^{-1} = \widehat{A}_L^{-1} \Delta_R \eqno(2.12) $$
using (2.8-10) and associativity of the algebra (2.7). 

The algebra (2.7) is associative under multiplication and the relations 
in (2.7) may be also expressed in a tensor product form 
$$ \widehat{R}(1) \widehat{A}_1 \widehat{A}_2 = - 
  \widehat{A}_2 \widehat{A}_1 \widehat{R}(1) \eqno(2.13) $$
where 
\begin{eqnarray*}
 \widehat{R}(x) 
  & = & (p + q^{-1}) \sum_i e_i^i \otimes e_i^i + 
        2x \sum_{i \neq j} (pq^{-1})^{i-1} e_i^i \otimes e_j^j + \\
  &   & (p - q^{-1}) \left ( \sum_{i > j} - \sum_{i < j} \right ) 
        e^i_j \otimes e^j_i 
\end{eqnarray*}
Here the elements of the matrix $e^k\hspace*{0.0001cm}_l$ are 
$$ \left (e^k\hspace*{0.0001cm}_l \right )^i_j = \delta^i_k \delta^j_l. 
                                                            \eqno(2.15)$$
The explicit form of $\widehat{R}(x)$ is 
$$ \widehat{R}(x) = \left ( \matrix{ 
 p + q^{-1} & 0 & 0 & 0 \cr
 0 & 2x & q^{-1} - p & 0 \cr
 0 & p - q^{-1} & 2xpq^{-1} & 0 \cr
 0 & 0 & 0 & p + q^{-1} \cr}
\right). \eqno(2.16) $$

In terms of the matrix elements Eq. (2.13) is of the form 
$$\widehat{R}^{ij} \hspace*{0.0001cm}_{{kl}} 
  \widehat{A}^k \hspace*{0.0001cm}_{m} 
  \widehat{A}^l \hspace*{0.0001cm}_{n} = - 
  \widehat{A}^i\hspace*{0.0001cm}_l \widehat{A}^i\hspace*{0.0001cm}_k 
  \widehat{R}^{kl}\hspace*{0.0001cm}_{ij}. \eqno(2.17) $$

Finally, we note that the algebra (2.7) and the $\widehat{R}$-matrix in 
(2.16) with $p = q$ and $x = - 1$ was given in Ref. 4 (Sec. III).

\noindent
{\bf III. TWO PARAMETER DEFORMATION OF THE GRASSMANN MATRIX SUPERGROUP} 

In this section, we consider 2x2 supermatrices whose diagonal elements are 
Grassmannian. We remark that the supergroup GL$(1\vert 1)$ whose nondiagonal 
elements are Grassmannian is $(p,q)$ deformed in Refs. 7 and 8. We will say 
that such supermatrices form the Grassmann supermatrix group and denote it by 
Gr$(1\vert 1)$. Explicity, a Grassmann 2x2 supermatrix $\widehat{\cal A}$ is 
of the form 
$$ \widehat{\cal A}  = \left ( \matrix{ \alpha & b \cr
                                       c & \delta \cr}
\right) \eqno(3.1) $$
with two odd (greek letters) and two even (latin letters) matrix elements. 
Even matrix elements commute with everything and odd matrix elements 
anti-commute among themselves. 

We begin with Manin's approach.$^3$ To do this, we consider the 
endomorphisms of a two-dimensional quantum superplane and 
its dual, denoted by $R_p[1|1]$ and $R_q^*[1|1]$, respectively. 
$$U = \left (\begin{array}{c} x \\ \xi \end{array} \right ) \in R_p[1|1] ~
      \Longleftrightarrow~ x \xi - p \xi x = 0,~ \xi^2 = 0, \eqno(3.2) $$
and its dual 
$$\widehat{U} = \left (\begin{array}{c} \eta \\ y \end{array} \right ) 
   \in R_q^*[1|1] ~\Longleftrightarrow~ 
   \eta^2 = 0,~ \eta y - q^{-1} y \eta = 0.    \eqno(3.3) $$

Suppose that the matrix elements of $\widehat{\cal A}$ (anti-)commute 
with the coordinates of $R_p[1|1]$ and $R_q^*[1|1]$. Then, the 
endomorphisms 
$$ \widehat{\cal A} : R_p[1|1] \longrightarrow R_q^*[1|1], \eqno(3.4a)$$
$$ \widehat{\cal A} : R_q^*[1|1] \longrightarrow R_p[1|1] \eqno(3.4b)$$
impose the following bilinear product relations among the generators 
of $\widehat{\cal A}$:
 $$ \alpha b = p^{-1} b \alpha, \qquad \alpha c = q^{-1}c \alpha, \eqno(3.5a)$$
 $$ \delta b = p^{-1}b \delta, \qquad \delta c = q^{-1}c \delta, \eqno(3.5b) $$
 $$\alpha \delta + \delta \alpha = 0, \qquad \alpha^2 = 0 = \delta^2,\eqno(3.5c)$$
 $$ bc = pq^{-1} cb + (p - q^{-1}) \delta \alpha \eqno(3.5d)$$
where $p$ and $q$ are non-zero complex numbers and $pq \pm 1 \neq 0$. 
These relations may be considered as a two parameter deformation of a Grassmann 
superalgebra on four elements $(\alpha, b, c, \delta)$ where $\alpha$ and 
$\delta$ are Grassmannian elements. This deformed algebra denoted by 
Gr$_{p,q}(1\vert 1)$. For $p = q$, one obtains the one parameter deformation 
of the generators of $\widehat{\cal A}$ that was given in Ref. 9. 

The inverse of $\widehat{\cal A}$ can be found as in Ref. 9 and it is of 
the form 
$$\widehat{\cal A}^{-1}  
 = \left(\matrix{- c^{-1} \delta b^{-1} & c^{-1} + 
                             c^{-1} \delta b^{-1} \alpha c^{-1} \cr \cr
                  b^{-1} + b^{-1} \alpha c^{-1} \delta b^{-1} & 
                               - b^{-1} \alpha c^{-1} \cr}
\right) \eqno(3.6) $$
provided that $b$ and $c$ are invertible. It is easy to verify that 
this is the proper right and left inverse of $\widehat{\cal A}$, i.e. 
$$\widehat{\cal A} \widehat{\cal A}^{-1} = I = 
 \widehat{\cal A}^{-1} \widehat{\cal A}. $$
Let 
$$ \widehat{\cal A}^{-1} = \left ( \matrix{ \alpha' & b' \cr
                                            c' & \delta' \cr}
\right). $$
Then, the matrix elements of $\widehat{\cal A}^{-1}$ satisfy the following 
relations
  $$ \alpha' b' = p b' \alpha', \qquad \alpha' c' = q c' \alpha', $$
  $$ \delta' b' = p b' \delta', \qquad \delta' c' = q c' \delta', \eqno(3.7)$$
  $$ \alpha' \delta' + \delta' \alpha' = 0 , \qquad 
     {\alpha'}^2 = 0 = {\delta'}^2, $$
  $$ b'c' = qp^{-1} c'b' + (q - p^{-1}) \delta' \alpha'. $$
Therefore, the matrix elements of $\widehat{\cal A}^{-1}$ satisfy the 
$(p^{-1},q^{-1})$-commutation relations while the matrix elements 
of $\widehat{\cal A}$ satisfy the $(p,q)$-commutation relations. 

The quantum (dual) superdeterminant of $\widehat{\cal A}$ is defined as 
$$s\widehat{D}_{p,q}(\widehat{\cal A}) = \widehat{D} 
    = c^{-1}b - c^{-1} \alpha c^{-1} \delta \\
    = pq^{-1}(bc^{-1} - \alpha c^{-1} \delta c^{-1}) \eqno(3.8)$$
which for $p = q$ is the same as $s\widehat{D}_q(\widehat{\cal A})$ 
in Ref. 9. The factor $pq^{-1}$ in (3.8) appeared because of the 
relation (3.5d). Note that the second equality in (3.8) is obtained 
by using the relation 
$$ bc^{-1} = qp^{-1} c^{-1}b - (q - p^{-1}) c^{-1} \delta \alpha c^{-1} 
    \eqno(3.9)$$
which in turn is obtained from equation (3.5d).

In general $\widehat{D}$, the quantum (dual) superdeterminant of 
$\widehat{\cal A}$, is not central but obeys the following commutation 
relations 
$$ \widehat{D} \alpha = pq^{-1} \alpha \widehat{D}, ~~
   \widehat{D} \delta = pq^{-1} \delta \widehat{D}, $$
$$ \widehat{D} b = pq^{-1} b \widehat{D}, ~~
   \widehat{D} c = pq^{-1} c \widehat{D}.   \eqno(3.10)$$

It is interesting that the quantum (dual) superdeterminant 
$\widehat{D}$ is not central while the quantum superdeterminant 
of a matrix $GL_{p,q}(1|1)$ is 7. Hovewer, it becomes 
central for $p = q$ as noted down in Ref. 9. 

Before passing to the next section, we remark that the interesting 
point in the construction of (3.6) is the fact that the dual 
superdeterminant $\widehat{D}$ is not necessarily central. 

\noindent
{\bf IV. THE $\widehat{R}$-MATRIX }

In this section, we give an $\widehat{R}$-matrix to obtain the 
relations (3.5). The algebra (3.5) is associative under multiplication 
and the relation (3.5) may be expressed in terms of a graded 
$\widehat{R}$-matrix condition, as with the quantum supermatrix. 
To this end, we use the tensoring convention 
$$ (\widehat{\cal A}_1)^{ij}\hspace*{0.0001cm}_{kl} 
   = (\widehat{\cal A} \otimes I)^{ij}\hspace*{0.0001cm}_{kl} 
   = (-1)^{k(j+l)} \widehat{\cal A}^i\hspace*{0.0001cm}_k 
         \delta^j\hspace*{0.0001cm}_l 
     = \widehat{\cal A}^i\hspace*{0.0001cm}_k \delta^j\hspace*{0.0001cm}_l,
\eqno(4.1\mbox{a})$$
$$ (\widehat{\cal A}_2)^{ij}\hspace*{0.0001cm}_{kl} 
   = (I \otimes \widehat{\cal A})^{ij}\hspace*{0.0001cm}_{kl} 
   = (-1)^{i(j+l)} \widehat{\cal A}^j\hspace*{0.0001cm}_l 
         \delta^i\hspace*{0.0001cm}_k. \eqno(4.1\mbox{a})$$

The explicit form of $\widehat{\cal A}_1$ and $\widehat{\cal A}_2$ is 
$$\widehat{\cal A}_1   
 = \left(\matrix{\alpha & 0 & b & 0 \cr 
                      0 & \alpha & 0 & b \cr
                      c & 0 & \delta & 0 \cr 
                      0 & c & 0 & \delta \cr }
\right), \eqno(4.2a)$$
$$\widehat{\cal A}_2 
 = \left(\matrix{- \alpha & - b & 0 & 0 \cr 
                      - c &- \delta & 0 & 0 \cr
                      0 & 0 & - \alpha & b \cr 
                      0 & 0 & c & - \delta \cr }
\right). \eqno(4.2b)$$
Then the associative algebra (3.5) is equivalent to equation 
$$ \widehat{R}(-1) \widehat{\cal A}_1 \widehat{\cal A}_2 
    = - \widehat{\cal A}_2 \widehat{\cal A}_1 \widehat{R}(-1) \eqno(4.3) $$
where
$$\widehat{R}(-1) 
 = \left(\matrix{ p + q^{-1} & 0 & 0 & 0 \cr 
                        0 & - 2  & q^{-1} - p & 0 \cr
                        0 & p - q^{-1} & - 2pq^{-1} & 0 \cr 
                        0 & 0 & 0 & p + q^{-1}  \cr }
\right). \eqno(4.4)$$
This $\widehat{R}(-1)$-matrix obtained from (2.14) with $x = - 1$. 
Here a 4x4 matrix in the form (4.2a) is labeled in the following way
$$M = \left(\matrix{ 
       M^{11}\hspace*{0.0001cm}_{11} & M^{11}\hspace*{0.0001cm}_{12} & 
       M^{11}\hspace*{0.0001cm}_{21} & M^{11}\hspace*{0.0001cm}_{22}   \cr 
       M^{12}\hspace*{0.0001cm}_{11} & M^{12}\hspace*{0.0001cm}_{12} & 
       M^{12}\hspace*{0.0001cm}_{21} & M^{12}\hspace*{0.0001cm}_{22}  \cr 
       M^{21}\hspace*{0.0001cm}_{11} & M^{21}\hspace*{0.0001cm}_{12} & 
       M^{21}\hspace*{0.0001cm}_{21} & M^{21}\hspace*{0.0001cm}_{22}  \cr 
       M^{22}\hspace*{0.0001cm}_{11} & M^{22}\hspace*{0.0001cm}_{12} & 
       M^{22}\hspace*{0.0001cm}_{21} & M^{22}\hspace*{0.0001cm}_{22}\cr }
\right) \eqno(4.5)$$
similar to Ref. 6. 

We have given the $(p,q)$-commutation relations which satisfied 
by the matrix elements of a Grassmannian matrix and a Grassmannian 
supermatrix, i.e. we made a two parameter deformation of the 
Grassmann matrix group Gr(2) and the supermatrix group Gr$(1|1)$. 
We obtained the Grassmannian quantum superdeterminant of a Grassmannian 
quantum supermatrix $(p,q)$-deformed case. Hovewer, it reduces to the 
case discussed in Ref. 9 for $p = q$. We have given an $\widehat{R}$-matrix 
which by use of a tensor product gives the $(p,q)$-commutation relations 
between the matrix elements of a Grassmannian supermatrix. 

\noindent
{\bf ACKNOWLEDGMENT}

This work was supported in part by TBTAK the Turkish Scientific and Technical 
Research Council. 

\noindent
{\bf APPENDIX: THE PROPERTIES OF THE $n$th POWER OF \\ GRASSMANN SUPERMARICES}

Here we will discuss the properties of the $n$th 
power of a Grassmann supermatrix. First we note that the product 
of two Grassmann supermatrices is not a Grassmann supermatrix, i.e., 
the matrix elements of a product $M = \widehat{M}\widehat{M}'$ 
do not satisfy (3.5). Hovewer, $\widehat{M}\widehat{M}' \in 
GL_{p,q}(1|1)$ if $\widehat{M}$ and $\widehat{M}'$ are two 
Grassmann supermatrices and $(b,c)$ ($(\alpha,\delta)$) pairwise 
commute (anti-commute) with $(b',c')$ ($(\alpha',\delta')$). 
So, we must consider the matrix elements of $\widehat{M}^n$ 
with respect to even and odd values of $n$. Let the $(2n-1)$-th 
power of $\widehat{M}$ be 
$$ \widehat{M}^{2n-1} = \left(\matrix {A_{2n-1} & B_{2n-1} \cr \cr
                                       C_{2n-1} & D_{2n-1}  \cr}
\right), ~n \geq 1.  \eqno(A1) $$
After some algebra, one obtains
$$ A_{2n-1} = \left \{<n>_{pq} \alpha + 
   p <n-1>_{pq} \delta \right \} (bc)^{n-1}, $$
$$ B_{2n-1} = \left \{bc + 
   p <n-1>_{p^2q^2} \alpha \delta \right \} (bc)^{n-2}b,  \eqno(A2)$$
$$ C_{2n-1} = \left \{cb + 
   q <n-1>_{p^2q^2} \delta \alpha \right \} (cb)^{n-2}c,  $$
$$ D_{2n-1} = \left \{<n>_{pq} \delta + 
   q <n-1>_{pq} \alpha \right \} (cb)^{n-1},  $$
where
$$ <N>_{pq} = \frac{1 - p^Nq^N}{1-pq}. \eqno(A3) $$
Now it is easy to show the following relations are satisfied:
$$ A_{2n-1} B_{2n-1} = p^{-(2n-1)}B_{2n-1}A_{2n-1} $$
$$ A_{2n-1} C_{2n-1} = p^{-(2n-1)}C_{2n-1}A_{2n-1} $$
$$ D_{2n-1} B_{2n-1} = q^{-(2n-1)}B_{2n-1}D_{2n-1} $$
$$ D_{2n-1} C_{2n-1} = q^{-(2n-1)}C_{2n-1}D_{2n-1}, \eqno(A4) $$
$$ A_{2n-1} D_{2n-1} + D_{2n-1} A_{2n-1} = 0,$$
$$A_{2n-1}^2 = 0 = D_{2n-1}^2,  $$
$$ B_{2n-1} C_{2n-1} =  p^{2n-1}q^{-(2n-1)}C_{2n-1}B_{2n-1} + 
            \left (p^{2n-1} - q^{-(2n-1)} \right ) A_{2n-1}D_{2n-1}. $$
Thus, $\widehat{M}^{2n-1}$ is a Grassmann supermatrix with deformation 
parameters $p^{2n-1}$ and $q^{2n-1}$, i.e. 
$\widehat{M}^{2n - 1} \in Gr_{p^{2n-1},q^{2n-1}}(1|1)$. 

Similarly, if we write the matrix $\widehat{M}^{2n}$, the (2n)-th 
power of $\widehat{M}$ as 
$$ \widehat{M}^{2n} = \left(\matrix {A_{2n} & B_{2n} \cr \cr
                                       C_{2n} & D_{2n}  \cr}
\right), ~n \geq 1  \eqno(A5) $$
where 
$$ A_{2n} = \left \{bc + p \frac{1 - pq}{1 + pq} <n>_{pq}<n-1>_{pq} 
              \alpha \delta \right \} (bc)^{n-1},$$
$$ B_{2n} = <n>_{pq} \left \{\alpha + p \delta \right \}(bc)^{n-1}b,\eqno(A6)$$
$$ C_{2n} = <n>_{pq} \{\delta + q \alpha\} (cb)^{n-1}c,  $$
$$ D_{2n} = \left \{cb + q \frac{1 - pq}{1 + pq} <n>_{pq}<n-1>_{pq} 
              \delta \alpha \right \} (cb)^{n-1}, $$
then the elements of $\widehat{M}^{2n}$ obey the following relations
$$ A_{2n} B_{2n} = p^{2n}B_{2n}A_{2n}, ~~
   A_{2n} C_{2n} = p^{2n}C_{2n}A_{2n} $$
$$ D_{2n} B_{2n} = q^{2n}B_{2n}D_{2n}, ~~
   D_{2n} C_{2n} = q^{2n}C_{2n}D_{2n},\eqno(A7) $$ 
$$ B_{2n} C_{2n} + p^nq^{-n}C_{2n} B_{2n} =0, $$
$$B_{2n}^2 = 0 = C_{2n}^2,  $$
$$ A_{2n} D_{2n} - D_{2n}A_{2n} = (p^{2n} - q^{-2n}) C_{2n}B_{2n}. $$
Thus the matrix $\widehat{M}^{2n}$ is a supermatrix in the form of 
$$ T = \left ( \matrix{     a & \beta \cr
                      \gamma  & d \cr}
\right) \eqno(A8) $$
with the deformation parameters $p^{2n}$ and $q^{2n}$. The $n$-th power 
of such a supermatrix and the relations between the matrix elements of 
$T^n$ can be found in Ref. 10 (Sec. 3).

\noindent
$^1$ L. D. Faddeev, N. Y. Reshetikhin and L. A. Takhtajan, Algebra Anal. 
     {\bf 1}, 178 (1987). \\
$^2$ J. Wess and B. Zumino, Nucl. Phys. B {\bf 70}, 39 (1974).\\
$^3$ Yu.I. Manin, Commun. Math. Phys. {\bf 123}, 163 (1989); \\
\hspace*{0.3cm}{\it Quantum groups and non-commutative geometry}, CRM-1561 (1988).\\
$^4$ E. Corrigan, B. Fairlie, P. Fletcher and R. Sasaki, 
  J. Math. Phys. {\bf 31}, 776 
\hspace*{0.3cm}{(1990).}\\
$^5$ W. Schmidke, S. Vokos and B. Zumino, Z. Phys. C. {\bf 48}, 249 (1990). \\
$^6$ J. Schwenk, B. Schmidke and S. Vokos, Z. Phys. C {\bf 46}, 643 (1990).\\
$^7$ L. Dabrowski and L. Wang, Phys. Lett. B {\bf 266}, 51 (1991).\\
$^8$ R. Chakrabarti and R. Jagannathan, J. Phys. A {\bf 24}, 5683 (1991).\\
$^9$ S. Celik and S. A. Celik, 
  Balkan Phys. Lett. {\bf 3(3)}, 188 (1995). \\
$^{10}$ S. Celik and S. A. Celik, Preprint MSUMB-95/03 (1995). 
(Balkan Phys. Lett. {\bf 1}, 32 (1997).)

\end{document}